\documentclass[a4paper]{article}
\usepackage{amssymb}
\usepackage{graphicx}

\newtheorem{theorem}{Theorem}[section]

\newtheorem{lemma}[theorem]{Lemma} 
 
\def\longequation{$$\vcenter\bgroup\advance\hsize by -9em%
\noindent\ignorespaces\refstepcounter{equation}}%
\makeatletter%
\def\endlongequation{\egroup\eqno(\theequation)$$\global\@ignoretrue}
\makeatother

\title{Graphs with no induced five-vertex \\ path or antipath}

\author{Maria Chudnovsky\thanks{Princeton University, Princeton, NJ 08544, USA
 E-mail: mchudnov@math.princeton.edu. Most of this work was conducted while the author was at Columbia University. Partially supported by NSF
grants DMS-1001091 and IIS-1117631.}
\and%
Louis~Esperet\thanks{CNRS, Laboratoire G-SCOP, University of Grenoble,
France. E-mail: louis.esperet@g-scop.grenoble-inp.fr}
\and%
Laetitia~Lemoine\thanks{Laboratoire G-SCOP, University of Grenoble,
France. E-mail: lae.lemoine@gmail.com}
\and%
Peter Maceli\thanks{Wesleyan University, Middletown CT 06459, USA. Most of this work was conducted while the author was at Columbia University. 
E-mail: pmaceli@wesleyan.edu.}
\and%
Fr\'ed\'eric~Maffray\thanks{CNRS, Laboratoire G-SCOP, University of
Grenoble, France.  E-mail: frederic.maffray@g-scop.grenoble-inp.fr.}
\and% 
Irena Penev \thanks{Department of Applied Mathematics and Computer
Science, Technical University of Denmark, Lyngby, Denmark.  Email:
ipen@dtu.dk.  A part of this work was conducted while the author was
at Universit\'e de Lyon, LIP, ENS de Lyon, Lyon, France.  Partially
supported by the LABEX MILYON (ANR-10-LABX-0070) of Universit\'e de
Lyon, within the program ``Investissements d'Avenir''
(ANR-11-IDEX-0007) operated by the French National Research Agency
(ANR), and by the ERC Advanced Grant GRACOL, project number 320812. 
\newline\newline
Authors~Esperet,~Maffray, Lemoine, and Penev were partially supported
by ANR project Stint under reference ANR-13-BS02-0007.}}

\begin{document}
\maketitle

\begin{abstract} 
\noindent We prove that a graph $G$ contains no induced $5$-vertex
path and no induced complement of a $5$-vertex path if and only if $G$
is obtained from $5$-cycles and split graphs by repeatedly applying
the following operations: substitution, split unification, and split
unification in the complement, where split unification is a new
class-preserving operation introduced here.
\end{abstract} 

\clearpage
\section{Introduction}

All graphs in this paper are finite and simple.  For fixed $n\ge 1$,
let $P_n$ denote the path on $n$ vertices, and for $n\ge 3$, let $C_n$
denote the cycle on $n$ vertices.  The graph $C_5$ is also called a
\emph{pentagon}.  The complement of a graph $G$ is denoted by
$\overline{G}$.  Given graphs $G$ and $F$, we say that $G$ is {\em
$F$-free} if no induced subgraph of $G$ is isomorphic to $F$.  Given a
family $\mathcal{F}$ of graphs, we say that a graph $G$ is {\em
$\mathcal{F}$-free} provided that $G$ is $F$-free for all $F \in
\mathcal{F}$.

A graph $G$ is \emph{perfect} if for every induced subgraph $H$ of $G$
the chromatic number of $H$ is equal to the maximum clique size in
$H$.  Chudnovsky, Robertson, Seymour, and Thomas \cite{CRST} solved the
long-standing and famous problem known as the Strong Perfect Graph
Conjecture by proving that a graph $G$ is perfect if and only if neither 
$G$ nor $\overline{G}$ contains an induced odd cycle of length at least 
five. 
 
There are various instances of the collection ${\cal F}$ such that 
${\cal F}$-free graphs are highly structured in a way that can be
described precisely; this fact is interesting in itself, and sometimes, 
it is also useful for solving various optimization problems on 
${\cal F}$-free graphs. The goal of this paper is to understand the structure of
$\{P_5,\overline{P_5}\}$-free graphs.  The motivation for this is
manifold: 

-- The class of $\{P_5, \overline{P_5}\}$-free graphs contains all
\emph{cographs} and all \emph{split graphs}.  Cographs are also known
as $P_4$-free graphs; their structure is very well understood (see for
example \cite{BLS,CPS}).  Split graphs are graphs whose vertex-set
can be partitioned into a clique and a stable set, and it is known
\cite{FolHam1,FolHam2} that they are exactly the $\{C_4,
\overline{C_4}, C_5\}$-free graphs.  

-- The class of $\{P_5, \overline{P_5}\}$-free graphs has already been
the object of much research.  Fouquet \cite{Fouquet} proved that the
study of this class can be reduced in a certain way (which we recall
in more detail below) to the study of $\{P_5, \overline{P_5},
C_5\}$-free graphs.  Moreover, it follows from the results of
Chv\'atal, Ho\`ang, Mahadev, and de~Werra~\cite{CHMW}, Giakoumakis and
Rusu \cite{GR}, and Ho\`ang and Lazzarato~\cite{HL} that several
optimization problems can be solved in polynomial time in the class of
$\{P_5, \overline{P_5}\}$-free graphs.  However, none of these results
gives (or attempts to give) a description of the structure of such
graphs.  

-- The class of $\{P_5, \overline{P_5}, C_5\}$-free graphs is a
subclass of the class of perfect graphs, and it is interesting to have a structure
theorem for this subclass since so far, no structure theorem has been 
proved for the class of all perfect graphs.

% \medskip

Before presenting our results, we need to introduce some notation and
definitions.  For a graph $G$, we denote by $V(G)$ its vertex-set and
by $E(G)$ its edge-set.  Given a set $S\subseteq V(G)$, let $N(S)$ be
the set of vertices in $V(G)\setminus S$ that have a neighbor in $S$.
Let $G[S]$ denote the subgraph of $G$ induced by $S$, and
let $G\setminus S$ denote the induced subgraph $G[V(G)\setminus S]$.
We say that a vertex $v$ in $V(G)\setminus S$ is \emph{complete} to
$S$ if $v$ is adjacent to every vertex of $S$, and that $v$ is
\emph{anticomplete} to $S$ if $v$ has no neighbor in $S$. A vertex of
$V(G)\setminus S$ that is neither complete nor anticomplete to $S$ is
\emph{mixed} on $S$.  Given two disjoint sets $S,T\subseteq V(G)$, we
say that $S$ is \emph{complete} to $T$ when every vertex of $S$ is
complete to $T$, and we say that $S$ is \emph{anticomplete} to $T$ when every
vertex of $S$ is anticomplete to $T$.

An \emph{anticomponent} of a set $S\subseteq V(G)$ is any subset of
$S$ that induces a component of the graph $\overline{G}[S]$.  A graph
$G$ is \emph{anticonnected} if $\overline{G}$ is connected.

A \emph{homogeneous set} is a non-empty set $S\subseteq V(G)$ such that every
vertex of $V(G)\setminus S$ is either complete or anticomplete to $S$.
A homogeneous set $S$ is \emph{proper} when $|S|\ge 2$ and $S\neq
V(G)$.  Let $G$ be a graph that admits a proper homogeneous set $S$,
and let $s$ be any vertex in $S$.  We can decompose $G$ into the two
graphs $G[S]$ and $G\setminus (S\setminus s)$. Since $S$ is a homogeneous set, 
we see that up to isomorphism, the latter graph is the same whatever the choice 
of $s$. Moreover, both $G[S]$ and
$G\setminus (S\setminus s)$ are induced subgraphs of $G$.   
The reverse operation, known as \emph{substitution}, can be defined as
follows.  Let $G$ and $H$ be two vertex-disjoint graphs and let $x$ be
a vertex in $G$.  Make a graph $G'$ with vertex-set $V(G\setminus
x)\cup V(H)$, taking the union of the two graphs $G\setminus x$ and
$H$ and adding all edges between $V(H)$ and the neighborhood of $x$ in
$G$.  Clearly, in $G'$, the set $V(H)$ is a homogeneous set,
$H=G'[V(H)]$, and $G$ is isomorphic to an induced subgraph of $G'$.
Moreover $V(H)$ is a proper homogeneous set if both $G$ and $H$ have
at least two vertices.  Thus, a graph $G$ is obtained by substitution
from smaller graphs if and only if $G$ contains a proper homogeneous
set.
A graph is \emph{prime} if it has no proper homogeneous set.

The following result about the structure of
$\{P_5,\overline{P_5}\}$-free graphs was proved by Fouquet in
\cite{Fouquet}.
\begin{theorem}[\cite{Fouquet}] \label{Fouquet} 
Every $\{P_5,\overline{P_5}\}$-free graph $G$ satisfies one of the
following properties:
\begin{itemize} 
\item $G$ contains a proper homogeneous set; 
\item $G$ is isomorphic to $C_5$; 
\item $G$ is $C_5$-free. 
\end{itemize} 
\end{theorem} 
Theorem~\ref{Fouquet} immediately implies that every $\{P_5,
\overline{P_5}, C_5\}$-free graph can be obtained by substitution
starting from $\{P_5,\overline{P_5},C_5\}$-free graphs and pentagons.
Furthermore, it is easy to check that every graph obtained by
substitution starting from $\{P_5,\overline{P_5},C_5\}$-free graphs
and pentagons is $\{P_5,\overline{P_5}\}$-free.  We remark that the
Strong Perfect Graph Theorem \cite{CRST} implies that a $\{P_5,
\overline{P_5}\}$-free graph is perfect if and only if it is
$C_5$-free.  Thus, every $\{P_5,\overline{P_5}\}$-free graph can be
obtained by substitution starting from $\{P_5,\overline{P_5}\}$-free
perfect graphs and pentagons.  In view of this, the bulk of this paper
focuses on prime $\{P_5,\overline{P_5},C_5\}$-free graphs
(equivalently: prime $\{P_5,\overline{P_5}\}$-free perfect graphs).

Our first result, Theorem~\ref{thm:struc}, states that every prime
$\{P_5,\overline{P_5},C_5\}$-free graph that is not split admits a
particular kind of partition.  Our second result,
Theorem~\ref{thm:decomp}, states that every prime
$\{P_5,\overline{P_5},C_5\}$-free graph that is not split admits a new
kind of decomposition, which we call a ``split divide'' (see section
\ref{sec:sdiv}).  Next, we reverse the split graph divide
decomposition and turn it into a composition that preserves the
property of being $\{P_5,\overline{P_5},C_5\}$-free.  We call this
composition ``split unification'' (see section \ref{sec:sunif}).
Finally, combining our results with Theorem~\ref{Fouquet}, we prove
that every $\{P_5,\overline{P_5}\}$-free graph is obtained by
repeatedly applying substitution, split graph unification, and split
graph unification in the complement starting from split graphs and
pentagons, and furthermore, we prove that every graph obtained in this 
way is $\{P_5,\overline{P_5}\}$-free (see Theorems~\ref{P5-P5c} and~\ref{cor}).

This paper results from the merging of the two (unpublished)
manuscripts \cite{CMP} and \cite{ELM} on the same subject; it combines
the proofs and results from these two manuscripts so as to present
them in the most succint way.

%%%%%%
%%%%%%
\section{Prime $\{P_5, \overline{P_5}, C_5\}$-free graphs}

Recall that a graph is \emph{split} if its vertex-set can be
partitioned into a stable set and a clique.  F\"oldes and Hammer
\cite{FolHam1,FolHam2} gave the following characterization of split
graphs (a short proof is given in~\cite[p.~151]{Gol}).
\begin{theorem}[\cite{FolHam1,FolHam2}]\label{thm:folham}
A graph is split if and only if it is $\{C_4, \overline{C_4},
C_5\}$-free.
\end{theorem}

\begin{lemma}\label{lem:abx}
In a $\{P_5, \overline{P_5}, C_5\}$-free graph $G$, let $A$ and $B$ be
non-empty and disjoint subsets of $V(G)$, and let $t$ be a vertex in
$V(G)\setminus(A\cup B)$ such that: 
\vspace{-.3cm}
\begin{quote}
$\bullet$ $t$ is anticomplete to $A$ and complete to $B$, \\
$\bullet$ every vertex in $B$ has a neighbor in $A$, and \\
$\bullet$ $A$ is connected. 
\end{quote}
\vspace{-.3cm}
Then some vertex of $A$ is complete to $B$.
\end{lemma}
\noindent{\it Proof.} Pick a vertex $a$ in $A$ with the maximum number
of neighbors in $B$.  Suppose that $a$ has a non-neighbor $y$ in $B$.
We know that $y$ has a neighbor $a'$ in $A$.  Since $A$ is connected,
there is a path $P=a_0$-$\cdots$-$a_k$ in $G[A]$ with $k\ge 1$,
$a_0=a'$ and $a_k=a$.  Choose $a'$ such that $k$ is minimal.  So $P$
is chordless and $y$ has no neighbor in $P\setminus \{a_0\}$.  Then
$k=1$, for otherwise $t, y, a_0, a_1, a_2$ induce a $P_5$.  By the
choice of $a$, since $y$ is adjacent to $a'$ and not to $a$, there is
a vertex $z$ in $B$ adjacent to $a$ and not to $a'$.  Then $a, z, t,
y, a'$ induce a $C_5$ or $\overline{P_5}$ (depending on the pair
$y,z$), a contradiction.  Thus $a$ is complete to~$B$.  $\Box$

\medskip

We say that a set, or a graph, is \emph{big} if it contains at least
two vertices.  
\begin{theorem}\label{thm:struc}
Let $G$ be a prime $\{P_5, \overline{P_5}, C_5\}$-free graph that
contains a $\overline{C_4}$.  Then there are pairwise disjoint subsets
$X_0, X_1, \ldots, X_m, Y_0, Y_1, \ldots, Y_m$, with $m\ge 2$, whose
union is equal to $V(G)$, such that the following properties hold,
where $X=X_0\cup X_1\cup\cdots\cup X_m$ and $Y=Y_0\cup
Y_1\cup\cdots\cup Y_m$:
\begin{enumerate}
\item\label{xy1}
For each $i\in\{1, \ldots, m\}$, $X_i$ is connected, $|X_i|\ge 2$,
$X_0$ is a (possibly empty) stable set, and $X_0, X_1, \ldots, X_m$
are pairwise anticomplete to each other.
\item\label{xy2}
For each $i\in\{1, \ldots, m\}$, $Y_i\neq\emptyset$, every vertex of
$Y_i$ is mixed on $X_i$ and complete to $X\setminus (X_i\cup X_0)$,
and $Y_0$ is complete to $X\setminus X_0$.
\item\label{xy3} 
$Y_0, Y_1, \ldots, Y_m$ are pairwise complete to each other.  (So each
anticomponent of $Y$ is included in some $Y_i$ with $i\in\{0,\ldots,
m\}$.)
\item\label{xy4}
No vertex of $X\setminus X_0$ is mixed on any anticomponent of $Y$.
\item\label{xy5}
For each $i\in\{1, \ldots, m\}$, $X_i$ contains a vertex that is
complete to $Y$.
\item\label{xy6}
Every vertex of $X_0$ is mixed on at most one anticomponent of $Y$.
\item\label{xy7}
For every big anticomponent $Z$ of $Y$, the set $X_Z$ of vertices of
$X_0$ that are mixed on $Z$ is not empty.  Moreover, if $Z$ and $Z'$
are any two distinct big anticomponents of $Y$, then $X_Z\cap X_{Z'}
=\emptyset$.
\item\label{xy8}
Each big anticomponent $Z$ of $Y$ contains a vertex that is
anticomplete to~$X_Z$.
\item\label{xy9}
If $Y$ is not a clique, there is a big anticomponent $Z$ of $Y$ such
that $X_Z$ is anticomplete to all big anticomponents of $Y\setminus
Z$.
\end{enumerate}
\end{theorem}
\noindent\emph{Proof.} Since $G$ contains a $\overline{C_4}$, there is
a subset $X$ of $V(G)$ such that $G[X]$ has at least two big
components.  We choose $X$ maximal with this property.  Let $X_1,
\ldots, X_m$ ($m\ge 2$) be the vertex-sets of the big components of
$G[X]$, and let $X_0=X\setminus (X_1\cup\cdots\cup X_m)$.  So
(\ref{xy1}) holds.  Let $Y=V(G)\setminus X$.  We claim that:

\begin{equation}\label{YY1}
\mbox{For every $y\in Y$ and $i\in\{1, \ldots, m\}$, $y$ has a 
neighbor in $X_i$.}
\end{equation}
Proof.  If $y$ has no neighbour in $X_i$, then $X\cup\{y\}$ induces a
subgraph of $G$ with at least two big components (one of which is
$X_i$), which contradicts the maximality of~$X$.  Thus (\ref{YY1})
holds.

\begin{longequation}\label{YY2}
For every vertex $y\in Y$, there is at most one integer $i$ in $\{1,
\ldots, m\}$ such that $y$ has a non-neighbor in $X_i$.
\end{longequation}
Proof.  Suppose that $y$ has a non-neighbor in two distinct components
$X_i$ and $X_j$ (with $1 \leq i,j \leq m$) of $X$.  For each $h\in\{i,j\}$, $y$ has a neighbor in
$X_h$ by (\ref{YY1}), and since $X_h$ is connected, there are adjacent
vertices $u_h, v_h\in X_h$ such that $y$ is adjacent to $u_h$ and not
to $v_h$.  Then $v_i, u_i, y, u_j, v_j$ induce a $P_5$, a
contradiction.  Thus (\ref{YY2}) holds.

An immediate consequence of Claims (\ref{YY1}) and (\ref{YY2}) is the
following.
\begin{longequation}\label{YY3}
For every vertex $y\in Y$, either $y$ is complete to $X\setminus X_0$,
or there is a unique integer $i\in\{1, \ldots, m\}$ such that $y$ is
complete to $X\setminus (X_i\cup X_0)$ and $y$ is mixed on $X_i$.
\end{longequation}

For each $i\in\{1, \ldots, m\}$, let $Y_i=\{y\in Y \mid y \mbox{ is
mixed on $X_i$}\}$, and let $Y_0=Y\setminus (Y_1\cup\cdots\cup Y_m)$.
By (\ref{YY3}), the sets $Y_0, Y_1, \ldots, Y_m$ are pairwise disjoint
and their union is $Y$.  For each $i\in\{1, \ldots, m\}$, since $G$ is
prime, $X_i$ is not a homogeneous set, so there exists a vertex in
$V(G)\setminus X_i$ that is mixed on $X_i$; by (\ref{xy1}), any such
vertex is in $Y$, and so $Y_i\neq\emptyset$.  Thus (\ref{xy2}) holds.

Now we prove (\ref{xy3}).  Suppose that $Y_i$ is not complete to $Y_j$
for some distinct $i,j\in\{0,\ldots,m\}$.  Let $y\in Y_i$ and $z\in
Y_j$ be non-adjacent.  Up to symmetry we may assume that $i\neq 0$,
say $i=1$.  Since $X_1$ is connected, there are adjacent vertices
$u_1$ and $v_1$ in $X_1$ such that $y$ is adjacent to $u_1$ and not to
$v_1$. By~(\ref{xy2}), $z$ is complete to $\{u_1, v_1\}$. Furthermore, 
by (\ref{xy2}), $Y_1$ is complete to $X_2$, and every vertex in $Y_j$ has a neighbor 
in $X_2$; thus, there exists a vertex $x_2 \in X_2$ such that $x_2$ is 
adjacent to both $y$ and $z$. By (\ref{xy1}), $x_2$ is non-adjacent to 
$u_1$ and $v_1$. But now $z, x_2, y, u_1, v_1$ induce a $\overline{P_5}$, 
a contradiction.  So the first sentence of (\ref{xy3}) holds.  
The second sentence is an immediate consequence of the first.  Thus 
(\ref{xy3}) holds.

Now we prove (\ref{xy4}).  Suppose on the contrary, and up to
symmetry, that a vertex $x$ in $X_1$ is mixed on some anticomponent
$Z$ of $Y$.  Since $Z$ is anticonnected, there are non-adjacent
vertices $y, z\in Z$ such that $x$ is adjacent to $y$ and not to $z$.
By~(\ref{xy2}), $z$ has a neighbor $u$ in $X_1$, so $z\in Y_1$.  Since
$X_1$ is connected, there is a path $u_0$-$\cdots$-$u_k$ in $G[X_1]$
with $u_0=u$, $u_k=x$ and $k\ge 1$.  Choose $u$ such that $k$ is
minimal.  By~(\ref{xy2}), $y$ has a neighbor $x_2$ in $X_2$, and since
$z\in Y_1$, $z$ is adjacent to $x_2$.  If $k=1$, then $x, y, z, u,
x_2$ induce a $C_5$ or $\overline{P_5}$ (depending on the pair $y,u$).
So $k\ge 2$.  The minimality of $k$ implies that $z$ is not adjacent
to $u_1$ or $u_2$, and $u$ is not adjacent to $u_2$.  Then $x_2, z, u,
u_1, u_2$ induce a $P_5$, a contradiction.

Now we prove (\ref{xy5}).  We observe that by (\ref{xy1}) and (\ref{xy2}), any vertex 
$t$ from a big component of $X\setminus X_i$ is complete to $Y_i$ and 
anticomplete to $X_i$, and so we can apply Lemma~\ref{lem:abx} to $X_i$, $Y_i$, and $t$.
It follows that some vertex $a$ of $X_i$ is complete to $Y_i$.
By~(\ref{xy2}), $X_i$ is complete to $Y\setminus Y_i$.  Thus $a$ is
complete to $Y$.

Now we prove (\ref{xy6}).  Suppose that a vertex $x$ in $X_0$ is mixed
on two anticompoments $Z_1$ and $Z_2$ of $Y$.  For each $j\in\{1,2\}$,
since $Z_j$ is anticonnected, there are non-adjacent vertices $y_j$
and $z_j$ in $Z_j$ such that $x$ is adjacent to $y_j$ and not to
$z_j$.  Then $y_1, z_1, x, z_2, y_2$ induce a $\overline{P_5}$, a
contradiction.

Now we prove (\ref{xy7}).  If $Z$ is any big anticomponent of $Y$,
then, since $G$ is prime, $Z$ is not a homogeneous set, and so there
exists a vertex of $V(G)\setminus Z$ that is mixed on $Z$.  The
definition of $Z$ and (\ref{xy4}) imply that any such vertex is in
$X_0$.  So $X_Z\neq\emptyset$.  The second sentence of (\ref{xy7})
follows directly from (\ref{xy6}).

Now we prove (\ref{xy8}).  Let $Z$ be a big anticomponent of $Y$.  By
(\ref{xy3}), $Z$ is included in one of $Y_0, Y_1, \ldots, Y_m$.
By~(\ref{xy2}) and (\ref{xy4}), some vertex $t$ of $X\setminus X_0$ is
complete to $Z$, and by (\ref{xy1}) $t$ is anticomplete to $X_Z$.
Hence we can apply Lemma~\ref{lem:abx} to $Z, X_Z$ and $t$ in the
complementary graph $\overline{G}$, and we obtain that some vertex in
$Z$ is complete (in $\overline{G})$ to $X_Z$.

Finally we prove (\ref{xy9}).  Suppose that $Y$ is not a clique, and
choose a big anticomponent $Z$ of $Y$ that minimizes the number of big
anticomponents of $Y$ that are not anticomplete to $X_Z$.  If this
number is $1$, then $Z$ satisfies the desired property.  So suppose
that this number is at least $2$, that is, there is a vertex $x\in
X_Z$ and a big anticomponent $Z'$ of $Y\setminus Z$ that contains a
neighbor of $x$.  There are non-adjacent vertices $y,z\in Z$ such that
$x$ is adjacent to $y$ and not to $z$.  By~(\ref{xy6}), $x$ is
complete to $Z'$.  Consider any $t\in X_{Z'}$; there are non-adjacent
vertices $y',z'\in Z'$ such that $t$ is adjacent to $y'$ and not to
$z'$.  If $t$ has any neighbor in $Z$, then, by~(\ref{xy6}), $t$ is
complete to $Z$, and then $z, x, t, z', y'$ induce a $\overline{P_5}$,
a contradiction.  Since this holds for any $t\in X_{Z'}$, we obtain
that $X_{Z'}$ is anticomplete to $Z$.  Now the choice of $Z$ implies
that there is a third big anticomponent $Z''$ of $Y$ (a big
anticomponent of $Y\setminus(Z\cup Z')$) such that some vertex $u$ of
$X_{Z'}$ has a neighbor $y''$ in $Z''$ and $X_Z$ is anticomplete to
$Z''$.  There are non-adjacent vertices $a,b\in Z'$ such that $u$ is
adjacent to $a$ and not to $b$.  Then $a, b, u, x, y''$ induce a
$\overline{P_5}$, a contradiction.  This completes the proof.  $\Box$

%%%

\section{The split divide}\label{sec:sdiv}

A \emph{split divide} of a graph $G$ is a partition $(A,B,C,L,T)$ of
$V(G)$ such that:
\begin{itemize}
\item
$|A|\ge 2$, $A$ is complete to $B$ and anticomplete to $C\cup T$, and
some vertex of $A$ is complete to $L$;
\item 
$L$ is a non-empty clique, every vertex of $L$ is mixed on $A$, and
$L$ is complete to $B\cup C$;
\item
$|C|\ge 2$, some vertex of $C$ is complete to $B$, and no vertex of
$C$ is mixed on any anticomponent of $B$;
\item
$T$ is a (possibly empty) stable set and is anticomplete to $C$.
\end{itemize}

\begin{figure}[htbp]
\begin{center}
\includegraphics[scale=1]{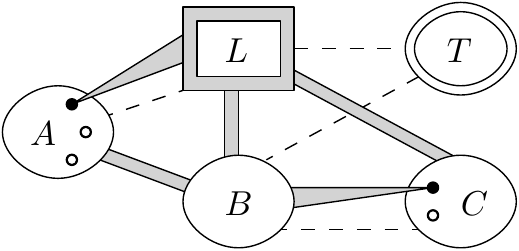}
\caption{A split divide. Adjacency between sets is as
  follows: gray means complete, no edge means anticomplete, and a
  dashed edge means arbitrary adjacency. Gray border around a set means 
  that the set is a clique, and white border means that the set is stable. 
\label{fig:divide}}
\end{center}
\end{figure}

Note that the sets $B$ and $T$ may be empty.  The split divide,
illustrated in Figure~\ref{fig:divide}, can be
thought of as a relaxation of the homogeneous set decomposition: a set
$X \subseteq V(G)$ is a homogeneous set in $G$ if no vertex in $V(G)
\setminus X$ is mixed on $X$; in the case of the split divide, the set
$A$ is not homogeneous, but all the vertices that are mixed on $A$ lie
in the clique $L$, and adjacency between $L$ and the rest of the graph
is heavily restricted.

%%%%%
\begin{theorem}\label{thm:decomp}
Let $G$ be a prime $\{P_5, \overline{P_5}, C_5\}$-free graph.  Then
either $G$ is a split graph or $G$ or $\overline{G}$ admits a split
divide.
\end{theorem}
\noindent{\it Proof.} \setcounter{equation}{0}
By Theorem~\ref{thm:folham} and up to
complementation, we may assume that $G$ contains a $\overline{C_4}$.
Consequently $G$ admits the structure described in
Theorem~\ref{thm:struc}, and we use it with the same notation.  All
items (i) to (ix) refer to Theorem~\ref{thm:struc}.

Suppose that $Y$ is a clique.  Let $A = X_1$, $L = Y_1$, $B = Y
\setminus Y_1$, $C = X_2\cup\cdots\cup X_m$ and $T = X_0$.  Then
$(A,B,C,L,T)$ is a split divide of $G$; this follows immediately from
the definition of the partition $X_0, X_1, \ldots, X_m, Y_0, Y_1,
\ldots, Y_m$, the fact that $Y$ is a clique, and
items~(\ref{xy1})--(\ref{xy5}).

Now suppose that $Y$ is not a clique.  We will show that
$\overline{G}$ admits a split divide.  By~(\ref{xy9}), we can choose a
big anticomponent $Z$ of $Y$ such that $X_Z$ is anticomplete to all
big anticomponents of $Y\setminus Z$.  By~(\ref{xy7}),
$X_Z\neq\emptyset$.  By~(\ref{xy3}), and up to relabeling, we may
assume that $Z\subseteq Y_0\cup Y_1$.  Hence $Z$ is complete to
$X_2\cup\cdots\cup X_m$, and every vertex of $X_1\cup (X_0\setminus
X_Z)$ is either complete or anticomplete to $Z$.  Let $K$ be the union
of all anticomponents of $Y$ of size~$1$.  So $K$ is a clique and is
complete to $Y\setminus K$.  Let:
\begin{eqnarray*}
A &=& Z;\\
L &=& X_Z;\\
B &=& \{x \in X_1\cup (X_0\setminus X_Z) \mid x \mbox{ is anticomplete
to $Z$}\};\\
C' &=& \{x \in X_1\cup (X_0\setminus X_Z) \mid x \mbox{ is complete to
$Z$}\}; \\
T &=& \{k\in K\mid k \mbox{ has a neighbor in $X_Z$}\}; \\
C &=& X_2\cup\cdots\cup X_m \cup (Y\setminus (Z\cup T))\cup C'.
\end{eqnarray*}
We claim that:
\begin{equation}\label{labc}
\mbox{$L$ is anticomplete to $B\cup C$.}
\end{equation}
Indeed, $X_Z$ ($=L$) is anticomplete to $X_1\cup\cdots\cup X_m$
because $X_Z\subseteq X_0$, and it is anticomplete to $X_0\setminus
X_Z$ because $X_0$ is a stable set.  Moreover, $X_Z$ is anticomplete
to every (big) anticomponent of $(Y\setminus K)\setminus Z$, by the
choice of $Z$, and it is anticomplete to $K\setminus T$ be the
definition of $T$.  Thus (\ref{labc}) holds.

\begin{equation}\label{cuv}
\mbox{No vertex of $C$ is mixed on any component of $B$.}
\end{equation}
For suppose that there is a vertex $c\in C$ and adjacent vertices
$u,v\in B$ such that $c$ is adjacent to $u$ and not to $v$.  Since
$X_0$ is a stable set and is anticomplete to $X_1$, we have $u,v\in
\{x \in X_1 \mid x \mbox{ is anticomplete to $Z$}\}$.  Since $c$ is
adjacent to $u$, we have $c\in (Y\setminus (Z\cup T)) \cup \{x \in X_1
\mid x \mbox{ is complete to $Z$}\}$.  Pick any $x\in X_Z$ and any
vertex $z\in Z$ adjacent to $x$.  By (\ref{labc}), $x$ is not adjacent to 
$c$.  Then $x, z, c, u, v$ induce
a $P_5$, a contradiction.  Thus (\ref{cuv}) holds.
 
\begin{equation}\label{kpc}
\mbox{$T$ is complete to $C$.}
\end{equation}
For suppose that there are non-adjacent vertices $t\in T$ and $c\in
C$.  Since $K$ is complete to $Y \setminus K$ and $T \subseteq K$, we
have that $c \notin Y \setminus (Z \cup T)$.  Thus, $c \in X_2
\cup\cdots\cup X_m \cup C'$.  By~(\ref{xy2}), $Y_0$ and $Y_1$ are
complete to $X_2 \cup\cdots\cup X_m$; since $Z \subseteq Y_0\cup Y_1$,
it follows that $Z$ is complete to $X_2 \cup\cdots\cup X_m$.  Thus,
$X_2 \cup\cdots\cup X_m\cup C'$ is complete to $Z$, and so $c$ is
complete to $Z$.  Further, since $X_2 \cup\cdots\cup X_m\cup C'
\subseteq X \setminus X_Z$ and $X_Z$ is anticomplete to $X \setminus
X_Z$ (because $X_Z \subseteq X_0$), we know that $c$ is anticomplete
to $X_Z$.  By the definition of $T$, $t$ has a neighbor $x$ in $X_Z$.
There are non-adjacent vertices $y,z\in Z$ such that $x$ is adjacent
to $y$ and not to $z$.  Since $t$ and $c$ are complete to $Z$, we see
that $t, c, y, z, x$ induce a $\overline{P_5}$, a contradiction.  Thus
(\ref{kpc}) holds.

Now we observe that:
\begin{itemize}
\item
$|A|\ge 2$ because $Z$ is big; $A$ is anticomplete to $B$ by the
definition of $B$; $A$ is complete to $C\cup T$ by~(\ref{xy2}); and
some vertex of $A$ is anticomplete to $L$ by~(\ref{xy8}).
\item 
$L$ is a non-empty stable set by (i) and (vii); every vertex of $L$ is
mixed on $A$ by the definition of $L$; and $L$ is anticomplete to
$B\cup C$ as shown in (\ref{labc}).
\item
$|C|\ge 2$ because $X_2\subseteq C$; some vertex of $C$ is
anticomplete to $B$ (every vertex of $X_2$ has this property); and no
vertex of $C$ is mixed on any component of $B$ as proved
in~(\ref{cuv}).
\item
$T$ is a clique and is complete to $C$ as proved in~(\ref{kpc}).
\end{itemize}
These observations mean that $(A, B, C, L, T)$ is a split divide in
$\overline{G}$.  This completes the proof.  $\Box$

\medskip

Let $G$ be a graph that admits a split divide $(A,B,C,L,T)$ as above,
let $a_0$ be a vertex of $A$ that is complete to $L$, and let $c_0$ be
a vertex of $C$ that is complete to $B$.  Let $G_1 = G[A\cup B \cup
\{c_0\} \cup L \cup T]$ and $G_2 = G[\{a_0\}\cup B \cup C \cup L \cup
T]$.  Then we consider that $G$ is decomposed into the two graphs
$G_1$ and $G_2$.  Note that $G_1$ and $G_2$ are induced subgraphs of
$G$ and each of them has strictly fewer vertices than $G$ since
$|A|\ge 2$ and $|C|\ge 2$.

\section{Split unification}\label{sec:sunif}

We can define a composition operation that ``reverses'' the split
divide decomposition.  Let $A,B,C,L,T$ be pairwise disjoint sets, and
assume that $A$ and $C$ are non-empty.  Let $a^*,c^*$ be distinct
vertices such that $a^*,c^* \notin A \cup B \cup C \cup L \cup T$.  \\
Let $G_1$ be a graph with vertex-set $A \cup B \cup L \cup T \cup
\{c^*\}$ and adjacency as follows:
\begin{itemize} 
\item 
$L$ is a (possibly empty) clique;
\item 
$T$ is a (possibly empty) stable set;
\item 
$A$ is complete to $B$ and anticomplete to $T$;
\item
Some vertex $a_0$ of $A$ is complete to $L$;
\item 
$c^*$ is complete to $B\cup L$ and anticomplete to $A \cup T$.
\end{itemize} 
Let $G_2$ be a graph with vertex-set $B \cup C \cup L \cup T \cup
\{a^*\}$ and adjacency as follows:
\begin{itemize} 
\item 
$G_2[B\cup L \cup T] = G_1[B\cup L \cup T]$;
\item 
$T$ is anticomplete to $C$; 
\item 
$L$ is complete to $B \cup C$; 
\item 
$a^*$ is complete to $B\cup L$ and anticomplete to $C \cup T$;
\item
Some vertex $c_0$ of $C$ is complete to $B$, and no vertex of $C$ is
mixed on any anticomponent of $B$.
\end{itemize} 

\begin{figure}[htbp]
\begin{center}
\includegraphics[scale=1]{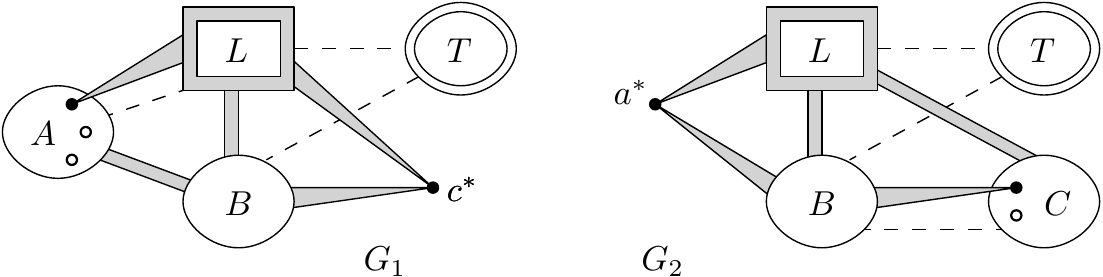}
\caption{A composable pair. \label{fig:unification}}
\end{center}
\end{figure}

Under these circumstances, we say that $(G_1,G_2)$ is a {\em
composable pair} (see Figure~\ref{fig:unification}).  The {\em split
unification} of a composable pair $(G_1,G_2)$ is the graph $G$ with
vertex-set $A \cup B \cup C \cup L \cup T$ such that:
\begin{itemize} 
\item 
$G[A \cup B \cup L \cup T] = G_1 \setminus c^*$; 
\item 
$G[B \cup C \cup L \cup T] = G_2 \setminus a^*$; 
\item 
$A$ is anticomplete to $C$ in $G$.
\end{itemize} 
Thus to obtain $G$ from $G_1$ and $G_2$, we ``glue'' $G_1$ and $G_2$
along their common induced subgraph $G_1[B\cup L \cup T]=G_2[B\cup L
\cup T]$, where $L\cup T$ induces a split graph (hence the name of the
operation).

We say that a graph $G$ is {\em obtained by split unification}
provided that there exists a composable pair $(G_1,G_2)$ such that $G$
is the split unification of $(G_1,G_2)$.  We say that $G$ is {\em
obtained by split unification in the complement} provided that
$\overline{G}$ is obtained by split unification.  We now prove that
every graph that admits a split divide is obtained by split
unification from smaller graphs.

\begin{theorem} \label{divuni} 
If a graph $G$ admits a split divide, then it is obtained from a
composable pair of smaller graphs (each of them isomorphic to an
induced subgraph of $G$) by split unification.
\end{theorem}
\noindent{\it Proof.} Let $G$ be a graph that admits a split divide.
Let $(A,B,C,L,T)$ be a split divide of $G$, let $a_0$ be a vertex of
$A$ that is complete to $L$, and let $c_0$ be a vertex of $C$ that is
complete to $B$.  Let $G_1=G[A\cup B\cup L \cup T\cup\{c_0\}]$.  Since
$|C| \geq 2$, we have $|V(G_1)|<|V(G)|$.  Let $G_2=G[B \cup C \cup L
\cup T\cup\{a_0\}]$.  Since $|A| \geq 2$, we have $|V(G_2)|<|V(G)|$.
Now $(G_1,G_2)$ is a composable pair, and $G$ is obtained from it by
split unification.  $\Box$

The split unification can be thought of as generalized substitution.
Indeed, we obtain the graph $G$ from $G_1$ and $G_2$ by first
substituting $G_1[A]$ for $a^*$ in $G_2$, and then reconstructing the
adjacency between $A$ and $L$ in $G$ using the adjacency between $A$
and $L$ in $G_1$.  We include $B$, $T$ and $c^*$ in $G_1$ in order to
ensure that split unification preserves the property of being
$\{P_5,\overline{P_5},C_5\}$-free.  In fact, we prove now something
stronger than this: split unification preserves the (individual)
properties of being $P_5$-free, $\overline{P_5}$-free, and $C_5$-free.

\begin{theorem}\label{compfree}
Let $(G_1,G_2)$ be a composable pair and let $G$ be the split
unification of $(G_1,G_2)$.  Then, for each $H\in\{P_5,
\overline{P_5}, C_5\}$, $G$ is $H$-free if and only if both $G_1$ and
$G_2$ are $H$-free.
\end{theorem}
\noindent{\it Proof.} We use the same notation as in the definition of
the split unification above.  First suppose that $G$ is $H$-free.
Observe that $G_1$ is isomorphic to the induced subgraph $G[A\cup
B\cup L\cup T\cup\{c_0\}]$, and $G_2$ is isomorphic to the induced
subgraph $G[B\cup C\cup L\cup T\cup\{a_0\}]$.  Hence $G_1$ and $G_2$
are $H$-free.  Now suppose that $G_1$ and $G_2$ are $H$-free and that
$G$ contains an induced copy of $H$.  Let $W$ be a five-vertex subset
of $V(G)$ such that $G[W] \simeq H$.  We claim that $W$ must contain
two non-adjacent vertices $b$ and $c$ with $b\in W\cap B$ and $c\in
W\cap C$.  For suppose the contrary.  Then $W\cap C$ is complete to
$W\cap (L\cup B)$ and anticomplete to $W\cap (A\cup T)$.  If $|W\cap
C|\ge 2$, then either $|W\cap C|\le 4$, so $W\cap C$ is a proper
homogeneous set in $G[W]$ (a contradiction since $H$ is prime), or
$W\subseteq C$, so $W$ is isomorphic to an induced subgraph of $G_2$
(a contradiction since $G_2$ is $H$-free).  So $|W\cap C|\le 1$, and
then $W$ is isomorphic to an induced subgraph of $G_1$ (where $c^*$
plays the role of the vertex in $W\cap C$ if there is such a vertex),
a contradiction since $G_1$ is $H$-free.  Therefore the claim holds.
By a similar argument, $W$ must contain two non-adjacent vertices $a$
and $\ell$ with $a\in W\cap A$ and $\ell\in W\cap L$.  Let $w$ be the
fifth vertex in $W$, so that $W=\{a,b,c,\ell,w\}$.  By the definition
of the split unification, $a,b,\ell,c$ induce a $P_4$ with edges $ab,
b\ell, \ell c$.  Consequently we must have one of the following two
cases: \\
(i) $W$ induces a $P_5$ or $C_5$.  So $w$ is anticomplete to $\{b,
\ell\}$ and has a neighbor in $\{a,c\}$.  Since $w$ is anticomplete to
$\{b,\ell\}$, it cannot be in $A, B, L$ or $C$, so it is in $T$.  But
then $w$ should be anticomplete to $\{a,c\}$.  \\
(ii) $W$ induces a $\overline{P_5}$.  So $w$ is adjacent to $a$ and
$c$ and has exactly one neighbor in $\{b,\ell\}$.  Since $w$ is
adjacent to $a$, it is not in $C\cup T$, and since it is adjacent to
$c$, it is not in $A$.  Moreover, since $w$ is adjacent to exactly one
of $b$ and $\ell$, it is not in $L$.  So $w\in B$, and so it is
adjacent to $\ell$ and, consequently, not to $b$.  Hence $b$ and $w$
lie in the same anticomponent of $B$, and $c$ is adjacent to exactly
one of them, a contradiction (to the last axiom in the definition of a
split unification).  $\Box$

\section{The main theorem} \label{section:main-thm} 

In this section, we use Theorem~\ref{Fouquet} and the results of the
preceding sections to prove Theorem~\ref{P5-P5c}, the main theorem of
this paper.

\begin{theorem}\label{P5-P5c} 
A graph $G$ is $\{P_5,\overline{P_5}\}$-free if and only if at least
one of the following holds:
\begin{itemize} 
\item 
$G$ is a split graph; 
\item 
$G$ is a pentagon; 
\item 
$G$ is obtained by substitution from smaller $\{P_5,
\overline{P_5}\}$-free graphs; 
\item 
$G$ or $\overline{G}$ is obtained by split unification from smaller
$\{P_5,\overline{P_5}\}$-free graphs.
\end{itemize} 
\end{theorem}
\noindent{\it Proof.} We first prove the ``if'' part.  If $G$ is a
split graph or a pentagon, then it is clear that $G$ is
$\{P_5,\overline{P_5}\}$-free.  Since both $P_5$ and $\overline{P_5}$
are prime, we know that the class of $\{P_5,\overline{P_5}\}$-free
graphs is closed under substitution, and consequently, any graph
obtained by substitution from smaller $\{P_5,\overline{P_5}\}$-free
graphs is $\{P_5,\overline{P_5}\}$-free.  Finally, if $G$ or
$\overline{G}$ is obtained by split unification from smaller
$\{P_5,\overline{P_5}\}$-free graphs, then the fact that $G$ is
$\{P_5,\overline{P_5}\}$-free follows from Theorem~\ref{compfree} and
from the fact that the complement of a $\{P_5,\overline{P_5}\}$-free
graph is again $\{P_5,\overline{P_5}\}$-free.
 
For the ``only if'' part, suppose that $G$ is a
$\{P_5,\overline{P_5}\}$-free graph.  We may assume that $G$ is prime,
for otherwise, $G$ is obtained by substitution from smaller
$\{P_5,\overline{P_5}\}$-free graphs, and we are done.  If some
induced subgraph of $G$ is isomorphic to the pentagon, then by
Theorem~\ref{Fouquet}, $G$ is a pentagon, and again we are done.  Thus
we may assume that $G$ is $\{P_5,\overline{P_5},C_5\}$-free.  By
Theorem~\ref{thm:decomp}, we know that either $G$ is a split graph,
or one of $G$ and $\overline{G}$ admits a split divide.  In the former
case, we are done.  In the latter case, Theorem~\ref{divuni} implies
that $G$ or $\overline{G}$ is the split unification of a composable
pair of smaller $\{P_5,\overline{P_5},C_5\}$-free graphs, and again we
are done.  $\Box$
 
As an immediate corollary of Theorem~\ref{P5-P5c}, we have the
following.
\begin{theorem}\label{cor}
A graph is $\{P_5,\overline{P_5}\}$-free if and only if it is obtained
from pentagons and split graphs by repeated substitutions, split
unifications, and split unifications in the complement.
\end{theorem}

Finally, a proof analogous to the proof of Theorem~\ref{P5-P5c} (but
without the use of Theorem~\ref{Fouquet}) yields the following result
for $\{P_5,\overline{P_5},C_5\}$-free graphs.  
\begin{theorem}\label{P5-P5c-C5} 
A graph $G$ is $\{P_5,\overline{P_5},C_5\}$-free if and only if at
least one of the following holds:
\begin{itemize}
\item 
$G$ is a split graph; 
\item 
$G$ is obtained by substitution from smaller $\{P_5, \overline{P_5},
C_5\}$-free graphs;
\item 
$G$ or $\overline{G}$ is obtained by split unification from smaller
$\{P_5,\overline{P_5},C_5\}$-free graphs.
\end{itemize} 
\end{theorem}

\section*{Acknowledgment} 

We would like to thank Ryan Hayward, James Nastos, Paul Seymour, and
Yori Zwols for many useful discussions.

\end{document}